\documentclass[5p,twocolumn]{elsarticle} 

\journal{Automatica}
\biboptions{numbers,sort&compress}

\usepackage{natbib}   

\usepackage{float}
\usepackage{caption}
\usepackage{subcaption}

\usepackage{xcolor}

\usepackage{amsfonts,latexsym,amsmath,amssymb}
\usepackage[mathscr]{eucal}
\usepackage{graphicx,color}
\usepackage{comment}
\usepackage{hyperref}
\usepackage[utf8]{inputenc}
\usepackage{mathtools}
\usepackage{hyperref}
\usepackage{comment}
\newtheorem{theorem}{Theorem}
\newtheorem{lemma}{Lemma}

\newtheorem{assumption}{Assumption}

\usepackage{lipsum}
\usepackage{mathtools}
\usepackage{cuted}
\usepackage{xifthen}

\newcommand{\overbar}[1]{\mkern 1.5mu\overline{\mkern-1.5mu#1\mkern-1.5mu}\mkern 1.5mu}

\catcode`\@=11
\def\downparenfill{$\m@th\braceld\leaders\vrule\hfill\bracerd$}
\def\overparen#1{\mathop{\vbox{\ialign{##\crcr\crcr
\noalign{\kern0.4ex}
\downparenfill\crcr\noalign{\kern0.4ex\nointerlineskip}
$\hfil\displaystyle{#1}\hfil$\crcr}}}\limits}
\catcode`\@=12

\def\NN{{\mathbb N}}    
\def\RR{{\mathbb R}}    

\DeclareMathOperator{\diag}{diag}
\DeclareMathOperator{\rank}{rank}

\def\He{\normalfont{\texttt{He}}}

\def\cO{{\mathcal{O}}}

\def\bx{{\textbf{x}}}

\def\cB{{\mathcal{B}}}

\def\cE{{\mathcal{E}}}
\newcommand\cV[3][]{\ifthenelse{\isempty{#1}}{\mathcal{V}_{#2} (#3)}{\mathcal{V}_{#2}^{#1} (#3)}}

\def\cN{{\mathcal{N}}}

\def\cL{{\mathcal{L}}}
\def\cG{{\mathcal{G}}}
\def\cV{{\mathcal{V}}}

\def\proof{%
    \par\vspace{0.410em}\noindent%
    \textbf{Proof.}\begingroup%
            }
\def\endproof{\null \ \null\hfill$\Box$\endgroup\\[0.410em]\null}
\usepackage{verbatim}

\begin{document}

\begin{frontmatter}

\title{Edges' Riemannian energy analysis for 
synchronization of 
multi-agent\\ nonlinear systems over undirected weighted graphs}

\tnotetext[t1]{Research partially funded by ANR Alligator (ANR-22-CE48-0009-01).}

\author[ECL]{A. Cellier-Devaux}
\ead{alexandre.cellier-devaux@ecl20.ec-lyon.fr}

\author[LAGEP]{D. Astolfi}
\ead{daniele.astolfi@univ-lyon1.fr}
\author[LAGEP]{V. Andrieu}
\ead{vincent.andrieu@univ-lyon1.fr}

\address[ECL]{Laboratoire Ampere Dpt. EEA of the Ecole Centrale de Lyon, Universit\'e
de Lyon, 69134 Ecully, France.}
\address[LAGEP]{Universit\'e Lyon 1, Villeurbanne, France -- CNRS, UMR 5007, LAGEPP, France. }

\begin{abstract}
In this note we investigate the problem of global exponential synchronization of multi-agent systems
described by non-linear input affine dynamics.
We consider 
the case of  networks described by
undirected connected graphs possibly without leader. 
We present a set of sufficient conditions
based on a Riemannian metric approach in order
to design a state-feedback distributed control law.
Then, we study the convergence properties of
the overall network.
By exploiting the properties
of the edge Laplacian we construct a  Lyapunov
function that allows to conclude 
global exponential synchronization of the overall network.
\end{abstract}


\end{frontmatter}

\section{Introduction}
\label{sec:intro}

The synchronization problem of  multi-agents systems has received increasing amounts of attention in the last decades. 
In this note, we focus on the case of 
homogeneous networks in which each agent is described by 
continuous-time nonlinear dynamics.
The main objective consists in  designing a
distributed feedback law  achieving synchronization of the network.
Such a feedback law must respect two constraints:
each agent can use only the information of its neighbours (i.e. the agents to which is communicating with), according to the given communication graph; 
the feedback law is zero on the synchronization manifold (i.e. when all the agents have the same behaviour).
Building up on the seminal works
in the linear context
\cite{scardovi2008synchronization}, 
notable results
have been 
made for nonlinear systems, 
see 
\cite{aminzare2014synchronization,isidori2014robust,pavlov2022nonlinear,zeng2016convergence,coraggio2018synchronization,giaccagli2021sufficient,giaccagli2024synchronization}
to cite a few.
Existing works can be organized from two different point of view:
the graph properties (e, g. directed or undirected, weighted or unweighted, with or without leader) or from the type of 
dynamics of each individual agent (e.g., linear or non-linear and which type of nonlinearities).
However, while for linear dynamics, the synchronization problem can be solved for many  
types of graph (see, e.g. \cite{li2009consensus,dutta2022strict}), for nonlinear dynamics a comprehensive result is still missing in the literature. 
To the best of our knowledge, the more general results
in terms of graph and agent's dynamics
are given by
\cite{li2012global} (directed strongly connected graphs)
and 
\cite{isidori2014robust} (directed weakly connected graphs).
However, both
results consider 
dynamics which are nonlinear in the state
but linear in the input. Furthermore, 
quadratic Lyapunov functions are employed.
Now, by re-organizing existing results
from the contraction point of view
\cite{simpson2014contraction,manchester2017control,andrieu2018some,tsukamoto2021contraction,davydov2022non},
it can be noted
that most of these works focused on the 
use of Euclidean metrics (or equivalently, 
quadratic Lyapunov functions)
or $L_1/L_\infty$ metrics (corresponding
to linear Lyapunov functions).
As a result, this severely restricts the individual agents' dynamics to be linear with respect to the input, with no coupling between the control input and the state
\cite{li2012global,aminzare2014synchronization,isidori2014robust,pavlov2022nonlinear,zeng2016convergence,coraggio2018synchronization}.

\color{black}

\color{black}
The first objective of this work is 
therefore to relax such a constrain 
by considering
generic input-affine nonlinear systems.
To this end we follow 
the recent results 
\cite{manchester2017control,giaccagli2021sufficient,giaccagli2024synchronization} where  a 
 Riemannian metric analysis is considered \cite{simpson2014contraction,manchester2017control}.
Such an analysis  
(also denoted in literature as control contraction metric
\cite{manchester2017control})
provide a set of sufficient
conditions for the design of an infinite-gain margin 
law \cite{giaccagli2023further} which is
essential to solve the synchronization problem
even for the linear case see \cite{li2009consensus,dutta2022strict}. 
The main limitation of \cite{giaccagli2021sufficient,giaccagli2024synchronization} is however that only 
undirected unweighted leader-connected graphs
are considered.
Therefore, 
the second objective of this work is to extend
such an analysis to more generic 
graph connections. In particular, 
we focus on the case of undirected
weighted graphs without the presence of a network leader.
We highlight that the study of leaderless network
is well motivated by several applications, such as 
electrical networks and 
grid-forming power-converters
\cite{dorfler2018electrical,colombino2019global}.
For this, we need to modify the  Lyapunov analysis
developed in \cite{giaccagli2021sufficient,giaccagli2024synchronization},
which is based on the use of  $N-1$ error dynamics
($N$ being the number of agents in the network)
with respect to the leader
and the properties of the Laplacian matrix
to conclude overall synchronization.
Such an approach, however, systematically 
fails when considering complex 
graphs which are not diagonalizable
(see, e.g. \cite{isidori2014robust}).
Furthermore, the quadratic Laypunov functions
proposed in 
\cite{li2012global,isidori2014robust}
cannot be directly employed because 
of absence of a strong linear stable component
in the agent's dynamics.

To overcome
such a difficulty, we revisit the 
graph theory and the 
properties of the incidence matrix connecting the Laplacian and the edge Laplacian in \cite[Chapter 2]{mesbahi2010graph}.
Such an edge Laplacian has been used, for instance, 
in \cite{zelazo2007agreement,dimarogonas2010stability,aminzare2014synchronization,zeng2016convergence, zeng2016edge} to solve the synchronization problem in the linear or Euclidean metric scenario. 
In particular, this approach focuses on the edges of the graph and a undirected network structure. 
The main idea consists in considering
a number of error dynamics  which is equivalent to the number of edges and not to the number of agents (e.g.,
$N-1$ in 
\cite{isidori2014robust,giaccagli2024synchronization}, 
or $N$  in \cite{dutta2022strict}).
To this end we use the properties of the edge Laplacian to create a specific edge matrix, and we use the control law in \cite[Theorem 1]{giaccagli2021sufficient}.
Then, we propose a new Lyapunov function
which corresponds to the 
Riemannian edge's energy.
We remark that such a Lypaunov
equation can be seen as  
an extension
(to the Riemannian and 
weighted cases)
of the 
Lyapunov function proposed in 
\cite{aminzare2014synchronization}.
Its derivative 
along the flow of solutions is always negative and
certify the overall network synchronization.

\section{Preliminaries and Problem Statement} \label{sec:problem}

A communication graph is described by a triplet $\cG=(\cV,\cE,A)$. With $\cV=\{v_1,v_2,...,v_N\}$ a set of $N \in \NN$ nodes (the $N$ agents), $\cE=\{(v_{k_1},v_{l_1}),...,(v_{k_m},v_{l_m})\} \subset \cV \times \cV$ a set of edges that models the interconnection between the nodes. For an undirected graph, if $(v_k,v_l)\in \cE$, then $(v_l,v_k)\in \cE$. 
We denote with $\cN_i$ the set of in-neighbors of node $i$, i.e. the set $\cN_i:=\{j\in\{1, \dots, N\} \, \lvert\, e_{ji}\in\cE\}$. 
Furthermore,  $A=(a_{ij}) \in \RR^{N \times N}$ is the adjacency matrix defined by the weight $a_{kl}\geq 0$ of each edge from node $k$ to node $l$. Notice that $a_{kl}=a_{lk}$ for an undirected graph. We denote by $L=(\ell_{ij}) \in \RR^{N \times N}$ the Laplacian matrix of the graph, defined as
$ \ell_{ij} = - a_{ij}$ for $i\ne j$, 
and $\ell_{ij} =  \sum_{k=1}^N a_{ik}$
for $i= j$,
For an undirected connected graph with $N$ nodes, $L$ is a positive semi-definite matrix and its spectrum is known and ordered as:
$$
0=\lambda_1<\lambda_2\leq ...\leq \lambda_M \text{, with } M=\rank(L)+1
$$
Given an undirected 
graph with $N$ nodes $(v_1,...,v_N)$ and $Q$ edges $(\varepsilon_1,...,\varepsilon_Q)$, by giving an arbitrary orientation for each edge of the graph, we denote by $E=(e_{ij}) \in \RR^{N \times Q}$ the incidence matrix defined as 
$$
e_{ij}=\left\{
\begin{array}{l}
-1 \text{ if }v_i \text{ is the initial node of edge } \varepsilon_j\\
+1 \text{ if }v_i \text{ is the terminal node of edge } \varepsilon_j\\
0 \text{ otherwize}
\end{array}
\right.
$$
Then, we denote $E_l$ the matrix with the positive term of $E$ and $E_k$ the one with the absolute value of the negative term of $E$, such as we have
$
E=E_l-E_k.
$
We denote by $w_i=a_{kl}>0$, for each edge $\varepsilon_i=(v_k,v_l)$, the weight of the edge. Thus, we introduce the weighting matrix $W$ defined as $W=\diag{(w_1,...,w_Q)}$. Notice that the Laplacian is directly related to the incidence matrix and the  weighting matrix as $L=EWE^\top$. Then, we introduce an edge variant of the Laplacian matrix, 
that is the edge Laplacian $L_e\in \RR^{Q \times Q}$ which is defined as 
\begin{equation}
    L_e=E^\top EW\,.
    \label{eqedge}
\end{equation}
By definition, the non-zero eigenvalues of the Laplacian matrix $L$ are eigenvalues of the edge Laplacian $L_e$.

We consider a homogeneous 
network of $N\geq 2$ agents
connected via a weighted undirected graph, 
in which the dynamics of each agent is described by the following nonlinear
ode
\begin{equation}
    \Dot{x}_i=f(x_i)+g(x_i)u_i, 
    \qquad i\in \{1, \ldots, N\},
    \label{dyn}
\end{equation}
where $x_i \in \RR^n$ is the state of the node $i$, $u_i \in \RR$ is the control action of node $i$, and $f:\RR^n\rightarrow\RR^n$ and $g:\RR^n \rightarrow\RR^n$ are $C^1$ functions which are the same for each agent. 
The main goal of this article 
is to design a distributed diffusive coupling
of the form 
 $$  
 u_i  = \sum_{j\in {\cN_i}}a_{ij}\big(\varphi(x_j)- \varphi(x_i)\big)
 $$
in order to guarantee global convergence of the distance between any two 
agents of the network, i.e. 
\begin{equation}
\noindent
 \sum  \|x_i(t)-x_j(t)\| 
 \leq  
 k \exp(-\mu t) \sum  \|x_i(0)-x_j(0)\| 
     \label{eq:global-sync}
\end{equation}
for any initial condition $x_i(0)\in \RR^n$, $i=\{1,\ldots, N\}$,
and for some $k,\mu>0$.
It's important to stress that 
for an agent $i$, the desired control law has to be designed taking into account only the state information of the agent itself and of the set of agents which are neighbours, 
namely $j\in \mathcal{N}_i$. Furthermore, 
$u_i=0$ whenever state agreement is achieved, 
i.e. when $x_i= x_j$ for any $j\in \mathcal{N}_i$.

\section{Properties of the edge Laplacian matrix} \label{sec:edge}

For an undirected  graph we discussed that we can give an arbitrary orientation for each edge of the graph. In this paper, we take a specific orientation for each edge to simplify the notation. In particular,  each edge between agents $k$ and $l$ (with $k<l$) is noted $(k,l)$. Then, we order the arcs of the form $(k,l)$ by increasing $k$ then by increasing $l$ for equal $k$. We then obtain the set of arcs of the network with $Q$ edges that we note $\Sigma=\{(k_1,l_1), (k_2,l_2), ..., (k_Q,l_Q)\}$. Thus, we defined the matrix $E$ with this orientation and we obtain the edge Laplacian $L_e$ with \eqref{eqedge}.
Then, the following property can be established.
Note that we don't require the graph to be connected.

\begin{lemma} \label{lemma:edge}
   Consider an undirected graph
   (not necessarily connected)
   with $N$ agents and $Q$ edges. Let $W$, resp. $E$, resp. $L=EWE^\top$, resp. $L_e=E^\top E W$,
   be the weight matrix, the incidence, Laplacian
   and edge Laplacian matrix. 
   Then   there exists a matrix $\Upsilon \in \RR^{Q \times Q}$ verifying 
    \begin{align}
       &E^\top L =\Upsilon E^\top \, ,
        \label{elue}
\\
        &W\Upsilon +\Upsilon^\top W \succ 0 \,.
        \label{wu}
    \end{align}
\end{lemma}

\proof\:
First of all, assume that $\ker(E)\neq \emptyset$ and assume its dimension is $\ell$. Let $(v_1, \dots, v_\ell)$ in $\RR^{Q}$ be an orthonormal basis of $\ker(E)$.
Let 
$$
\Upsilon = E^\top EW+\mu \sum_{i=1}^\ell v_iv_i^\top
$$ 
where $\mu$ is a  real number that will be selected later on.
Note that 
\begin{align*}
\Upsilon E^\top &= E^\top EWE^\top +\mu\sum_{i=1}^\ell v_iv_i^\top E^\top
\\
&= E^\top L + \mu\sum_{i=1}^\ell v_i(Ev_i)^\top =E^\top L.
\end{align*}
Moreover, 
$
    W\Upsilon + \Upsilon^\top W = Q_1 + \mu Q_2
$
where 
$$
Q_1=2WE^\top EW,\qquad Q_2 =  W\sum_{i=1}^\ell v_iv_i^\top
    +\sum_{i=1}^\ell v_iv_i^\top W.
$$
Let $w$ in $\RR^Q$ be such that $w^\top Q_1w=0$ and $w\neq 0$. Hence  $EWw=0$ and this implies 
$Ww = \sum_{j=1}^\ell a_j v_j$ for $\ell$ real numbers $a_i$ not all identically zero.
Consequently, since $W$ is diagonal,
\begin{align*}
    w^\top Q_2w 
       &=  2 (Ww)^\top \sum_{i=1}^\ell v_iv_i^\top W^{-1} Ww,\\
    &=  2 (\sum_{j=1}^\ell a_j v_j^\top)\sum_{i=1}^\ell v_iv_i^\top W^{-1} (\sum_{j=1}^\ell a_j v_j),\\
    &= 2(\sum_{j=1}^\ell a_j v_j^\top)W^{-1}(\sum_{j=1}^\ell a_j v_j)>0,
\end{align*}
where the last inequality is obtained since $W$ is composed of positive element.
In conclusion, we have the implication 
$$(w^\top Q_1 w =0, w\neq 0) \ \Rightarrow\  w^\top Q_2 w>0.$$ Hence, with Finsler's lemma, there exists $\mu$ such that 
$
    W\Upsilon + \Upsilon^\top W \succ 0
$
concluding the proof.
\endproof{}

\section{Main Result} \label{sec:results}

In this section we consider a network represented by a graph $\mathcal G$
and satisfying the next assumption
of $N$ identical agents described by the
nonlinear dynamics \eqref{dyn}.  

\begin{assumption}\label{ass:graph}
The graph $\mathcal G$ is undirected and connected.
\end{assumption}

The objective is to design a distributed control law solving
the synchronization problem as in 
\eqref{eq:global-sync}.
We recall the definition 
of Lie derivative of  
$C^1$
tensor 
$P : \RR^n \to \RR^{n \times n}$ along a 
$C^1$
vector field $f:\RR^n \to \RR^n$
 defined as
$$
\cL_f P(x) = \mathfrak{d}_fP(x)+P(x) \tfrac{\partial f}{\partial x}(x)+\tfrac{\partial f}{\partial x}(x)^\top P(x)
$$
where
$
\mathfrak{d}_fP(x)=\lim_{h \to 0} \tfrac{P(x+hf(x))-P(x)}{h}
$
for all $x \in \RR^n$.
Then, following
\cite{giaccagli2021sufficient,giaccagli2024synchronization}, we introduce 
the following set of  assumptions.

\begin{assumption}\label{ass:riem}
There exists a $C^1$ function $P : \RR^n \rightarrow \RR^{n \times n}$ taking symmetric positive definite
values and 
real numbers $\mu,\rho,\underline{p},\overline{p}>0$
satisfying for all $x \in \RR^n$
     \begin{equation}
        \begin{aligned}
            & \cL_f P(x) - \rho \,P(x)g(x)g(x)^\top P(x)\preceq- 2\mu P(x),\\
            & \underline{p}I_n \preceq P(x) \preceq \overline{p}I_n.
            \label{a}
        \end{aligned}
        \end{equation}
\end{assumption}

\begin{assumption}\label{ass:killing}
The following Killing vector property 
and integrability conditions are verified
\begin{equation}\label{bc}
    \cL_g P(x)=0,
    \quad
       \cfrac{\partial \alpha}{\partial x}(x)=g(x)^\top P(x), \quad \forall x \in \RR^n,
\end{equation}
for some $C^2$ function 
$\alpha:\RR^n \to \RR$.
\end{assumption}

As detailed in \cite{giaccagli2023further},
Assumptions
~\ref{ass:riem} and
~\ref{ass:killing}
are sufficient to establish 
the existence of an infinite-gain margin  feedback 
law  making contractive a system of the form 
\eqref{dyn}.
Further comments are given at the end of the section. 
We can state now the following theorem.

\begin{theorem}\label{the1}
 Consider a network described by $N$ agents \eqref{dyn}
 and suppose that 
 Assumptions~\ref{ass:graph}-\ref{ass:killing} hold.
 Assume moreover, that the functions $f$, $\alpha$ are globally Lipschitz and $g$ is globally bounded.
 Let $\beta^\star = \frac{\rho\overline{w}}{2\underline{\lambda}}$
 with $\rho$ satisfying \eqref{a}, $\overline{w}$ being the maximal eigenvalue of the matrix $W$, and $\underline{\lambda}$ the minimal eigenvalue of the matrix $W\Upsilon$ with $\Upsilon$ satisfying \eqref{elue}.
Then, for any  $\beta\geq\beta^\star$, 
the distributed feedback
    \begin{equation}
        u_i=-\beta \sum_{j=1}^N \ell_{ij} \alpha(x_j)
     = \beta \sum_{j\in {\cN_i}}a_{ij}\Big(\alpha(x_j)- \alpha(x_i)\Big)   
        \label{ui}
    \end{equation}
with $\alpha$ satisfying \eqref{bc}, 
solve the global synchronization problem
for the network \eqref{dyn}.
More precisely, \eqref{eq:global-sync} holds
with $\mu$ given by Assumption~\ref{ass:riem}.
\end{theorem}

\proof \;
Let $\Sigma=\{(k_1,l_1), (k_2,l_2), ..., (k_Q,l_Q)\}$ be the set of arcs of the network  as defined in Section \ref{sec:edge}. 
For a given initial condition  $(x_1(0),\dots,x_{N}(0)$ in $\RR^{Nn}$ let $t\mapsto(x_1(t),\dots,x_{N}(t))$ be the solution of \eqref{dyn}. This one is defined for all positive time since by assumption, the nonlinear system is globally Lipschitz.
For any $i \in \{1,\ldots, Q\}$, 
we define 
the $C^2$ function $(s,t)\in[0,1]\times\RR_+\mapsto\Gamma_i(s,t)\in\RR^n$ solution of
\begin{align}
    \dfrac{\partial \Gamma_i}{\partial t}(s,t)&=f(\Gamma_i(s,t))+g(\Gamma_i(s,t))\widetilde{u}_i(s,t),
    \label{gamma}
\\
    \widetilde{u}_i(s,t)&=-\beta \left[\sum_{j=1}^Q  \Upsilon_{ij}\alpha(\Gamma_{j}(s,t))+ \sum_{l=1}^N \Omega_{il}\alpha(x_l(t)) \right],
    \label{utilde}
\end{align}
and with $\Gamma_i(0,s)=(1-s)x_{k_i} + sx_{l_i}$,
where $\Upsilon=(\Upsilon_{ij})$ is solution of \eqref{elue} and $\Omega=(\Omega_{il}) \in \RR^{Q \times N}$ is defined by $\Omega=\frac{1}{2}(\lvert E^\top \lvert L - \Upsilon \lvert E^\top \lvert)$, 
with the notation $\lvert A\rvert:=(\lvert a_{ij}\rvert)$, for any matrix 
$A=(a_{ij})$.
The function $\Gamma_i$ is defined as any $C^2$
path between the points $x_{k_i}$
and $x_{\ell_i}$ satisfying the dynamics 
\eqref{gamma},
  consistently
with the agent dynamics \eqref{dyn}.
Consider now the  $C^1$
function
$V:\RR_+\mapsto \RR_+$ defined by
\begin{equation}
    \label{eq:Lyap_function}
\begin{aligned}
V(t) & :=\sum_{i=1}^Q w_i
V_i(t)\;,\\
V_i(t) & := 
\int_{0}^{1} \left( \cfrac{\partial \Gamma_i^\top}{\partial s}(s,t) P(\Gamma_i(s,t)) \cfrac{\partial \Gamma_i}{\partial s}(s,t) \right) \,ds ,
\end{aligned}
\end{equation}
where $P$ is such that \eqref{a} and \eqref{bc} are verified. 
It corresponds to the sum of the energies
$V_i$ associated to the paths $\Gamma_i$
between any two agents connected via an edge $i\in \{1,\ldots, Q\}$.
Note that using the second inequality in 
\eqref{a} and the fact that for all $t\geq 0$,
$\Gamma_{i}(0,t)=x_{k_i}(t)$
and $\Gamma_{i}(1,t)=x_{l_i}(t)$, we obtain
$$
\begin{aligned}
    V_i(t) & \geq \underline p
\int_{0}^{1} \cfrac{\partial \Gamma_i^\top}{\partial s}(s,t) \cfrac{\partial \Gamma_i}{\partial s}(s,t)  \,ds
\geq 
\underline p\, \Vert x_{k_i}(t) - x_{l_i}(t) \Vert^2.
\end{aligned}
$$
So for all $t\geq 0$ it yields
\begin{align*}
    V(t) &  \geq \underline{p}\sum_{i=1}^Q w_i\, \Vert x_{k_i} (t)- x_{l_i} (t)\Vert^2 \\
    & = \underline{p} \bx(t)^\top (EW \otimes I_n)(E^\top \otimes I_n)\bx(t) = \underline{p} \bx^\top (L\otimes I_n) \bx(t)
\end{align*}
with $\bx=(x_1^\top, ..., x_N^\top)^\top$. 
The graph being connected, it implies  $\ker(L)=\text{Vect}(1_n)$. 
Hence, 
$
\ker(L\otimes I_n)=\text{Vect}(1_n\otimes I_n)
$
Consequently there exists a positive real number $\underline \alpha$ such that,  
\begin{equation}\label{UpperBound}
V(t)\geq \underline\alpha \sum_{i,j\in\{1, \ldots, N\}} \|x_i(t)-x_j(t)\|^2 
,\quad  \forall t\geq 0.
\end{equation}
On the other hand, the time-derivative of $V$ satisfies,
\begin{align*}
    \Dot{V}(t)=\sum_{i=1}^Q w_i \int_{0}^{1}
    & \biggl( \cfrac{\partial^2 \Gamma_i^\top}{\partial s \partial t}(s,t)P(\Gamma_i(s,t))\cfrac{\partial \Gamma_i}{\partial s}(s,t)\\
    & + \cfrac{\partial \Gamma_i^\top}{\partial s}(s,t)\mathfrak{d}_{f+g\widetilde{u}_i}P(\Gamma_i)\cfrac{\partial \Gamma_i}{\partial s}(s,t) \\
    & + \cfrac{\partial \Gamma_i^\top}{\partial s}(s,t)P(\Gamma_i(s,t))\cfrac{\partial^2 \Gamma_i}{\partial s \partial t}(s,t) \biggr) \,ds. 
\end{align*}
Using the definition of $\Gamma$
in \eqref{gamma}
and $\tilde u$ in
\eqref{utilde}, 
the term $\partial^2\Gamma / \partial s\partial t$
is computed as 
\begin{align*}
    \cfrac{\partial^2 \Gamma_i}{\partial s \partial t}(s,t)  = &\cfrac{\partial}{\partial s} \left[ f(\Gamma_i(s,t))+ g(\Gamma_i(s,t))\widetilde{u}_i \right]\\
     =&\cfrac{\partial f(\Gamma_i(s,t))}{\partial x}\cfrac{\partial \Gamma_i}{\partial s}(s,t) + \cfrac{\partial g(\Gamma_i(s,t))}{\partial x}\cfrac{\partial \Gamma_i}{\partial s}(s,t)\widetilde{u}_i\\
    & -\beta g(\Gamma_i(s,t))\cfrac{\partial \alpha(x)}{\partial x}\left( \sum_{j=1}^Q \Upsilon_{ij}\cfrac{\partial \Gamma_j}{\partial s}(s,t) \right)\\
     =& \left( \cfrac{\partial f(\Gamma_i(s,t))}{\partial x} + \cfrac{\partial g(\Gamma_i(s,t))}{\partial x}\widetilde{u}_i \right)\cfrac{\partial \Gamma_i}{\partial s}(s,t) \\
    & - \beta g(\Gamma_i)g(\Gamma_j)^\top P(\Gamma_j)\left( \sum_{j=1}^Q \Upsilon_{ij}\cfrac{\partial \Gamma_j}{\partial s}(s,t) \right).
\end{align*}
Then, combining  \eqref{bc} with 
the previous expression in the derivative of $V$
yields for all $t\geq0$
\begin{align*}
    & \Dot{V}(t) = \sum_{i=1}^Q w_i \int_{0}^{1} \biggl(\cfrac{\partial \Gamma_i^\top}{\partial s} \biggl[ \cL_f P(\Gamma_i) + \widetilde{u}_i \cL_g P(\Gamma_i)
     \biggr] \cfrac{\partial \Gamma_i}{\partial s} \\
    & - 2\beta \sum_{j=1}^Q \Upsilon_{ij} \text{He}\left[\cfrac{\partial \Gamma_i^\top}{\partial s}P(\Gamma_i)g(\Gamma_i)g(\Gamma_j)^\top P(\Gamma_j)\cfrac{\partial \Gamma_j}{\partial s}\right] \biggr) \,ds,
\end{align*}
with the compact notation 
$\He(A) = \tfrac{A+ A^\top}{2}$ for any element $A$.
By introducing the   notation
$ \Bar{\Lambda}(s,t)=(\Lambda_i^\top,...,\Lambda_Q^\top)^\top$
and  $  \Lambda_i(s,t)=g(\Gamma_i(s,t))^\top P(\Gamma_i(s,t)) \cfrac{\partial \Gamma_i}{\partial s}(s,t),$
 we obtain
\begin{align*}
    \Dot{V}(t) \leq &- 2\mu V(t) + k \sum_{i=1}^Q \int_{0}^{1} w_i\Lambda_i^\top \Lambda_i \, ds\\
    & - \beta \sum_{i=1}^Q \sum_{j=1}^Q \int_{0}^{1} w_i\Upsilon_{ij}[\Lambda_i^\top\Lambda_j + \Lambda_j^\top \Lambda_i]\,ds ,
\end{align*}
with \eqref{a} and \eqref{bc}, and therefore
\begin{align*}
    \Dot{V}(t) &\leq -2 \mu V(t) - \int_{0}^{1} \Bar{\Lambda}^\top\left[\beta (W\Upsilon+\Upsilon^\top W)-kW\right]\Bar{\Lambda}\,ds  \\
    & \leq -2\mu V(t) ,
\end{align*}
for any $\beta \geq \beta^\star$, 
$\beta^\star=\tfrac{\rho\overbar{w}}{2\underline{\lambda}}$.
This shows exponential convergence of $V$ to $0$. With \eqref{UpperBound}, it ensures the synchronization of the agents.
\endproof{}

\noindent
\textbf{Remark}.
Note that in the proof
of Theorem~\ref{the1}, 
    the function $\widetilde{u}_i$ in \eqref{utilde} is designed by the equations verified by $\Gamma_i$ for $s=0$ and $s=1$ for all
    $i\in \{1, \ldots, Q\}$. Indeed these equations yield $\widetilde{u}_i(0,t) =u_{k_i}(t)$ and $\widetilde{u}_i(1,t)=u_{l_i}(t)$. And with \eqref{ui}, it leads to
      \begin{align*}
        \widetilde{u}_i(0,t) &=  \sum_{m=1}^N (\Upsilon E_k^\top+\Omega)_{i,m} \alpha(x_m) \\
               u_{k_i}(t) &= \sum_{m=1}^N \ell_{k_i,m}\alpha(x_m) = \sum_{m=1}^N (E_k^\top L)_{i,m}\alpha(x_m) 
    \end{align*}
    so that $
           \widetilde{u}_i(0,t)=u_{k_i}(t) \Leftrightarrow {E_k^\top L=\Upsilon E_k^\top +\Omega}$
    and with the same equation for $s=1$, we obtain
    $        \widetilde{u}_i(1,t)=u_{l_i}(t) \Leftrightarrow  {E_l^\top L=\Upsilon E_l^\top +\Omega}.$
Subtracting 
the last two equations we recover
 \eqref{elue}.

\medskip
\noindent
\textbf{Remark.}
We highlight that the controller \eqref{ui}
is a diffusive coupling feedback which compares the difference 
between the $i$-th agent state and its $j$-th neighbours. 
As already remarked in 
\cite{giaccagli2021sufficient,giaccagli2023further,giaccagli2024synchronization}, 
in the case of linear dynamics of the form  
$$
\dot x_i = Ax_i+ Bu_i,
$$
Assumptions~\ref{ass:riem} -\ref{ass:killing}
boil down to the stabilizability of the pair $A,B$, and 
the matrix function $P$ is the solution to the Algebraic Riccati Inequality
$$
PA +A^\top P - \rho PBB^\top P \preceq - 2 \mu P.
$$
As a consequence, the function $\alpha$ can be selected as
$\alpha(x) = Kx$, with $K = B^\top P$,
and the distributed feedback law \eqref{ui}
simply reads
$u_i = -\beta \sum_{j=1}^N \ell_{ij} B^\top P x_j$, 
thus recovering 
standard existing solution for the linear context
\cite{scardovi2008synchronization,li2009consensus,dutta2022strict}.
In this linear case, 
the Lyapunov function certifying 
global expeontnial synchronization is
given by 
$$
V = e^\top (W \otimes P)e={\bx}^\top(EWE^\top \otimes P){\bx}
$$ 
where $e:=(E^\top \otimes I_n) {x}\in \RR^{Qn}$ with ${{\bx}}=(x_1, \ldots, x_N)$ being the vector of the state of the $N$ agents.
\color{black}

\section{Illustration} \label{sec:test}

The conditions of Assumption~\ref{ass:riem} and
~\ref{ass:killing} are based on contraction theory.
As a consequence, 
although they are very hard to verify analytically, 
they provide a feedback design which is robust by construction.
For this reasons, many recent works studied how to approximate
contraction-based controllers (also denoted in literature
as control contraction matrices \cite{manchester2017control}).
See, for instance, 
\cite{tsukamoto2021contraction}.
In particular, as already shown in 
\cite{giaccagli2024synchronization}, 
the integrability and Killing vector conditions
can be relaxed at the price of obtaining a 
feedback design which allows to achieve
only semi-global state synchronization
(where the domain of attraction has to be intended
with respect to a asymptotic behaviour of the 
overall network).

\begin{figure}    \centering
    \includegraphics[scale=0.1]{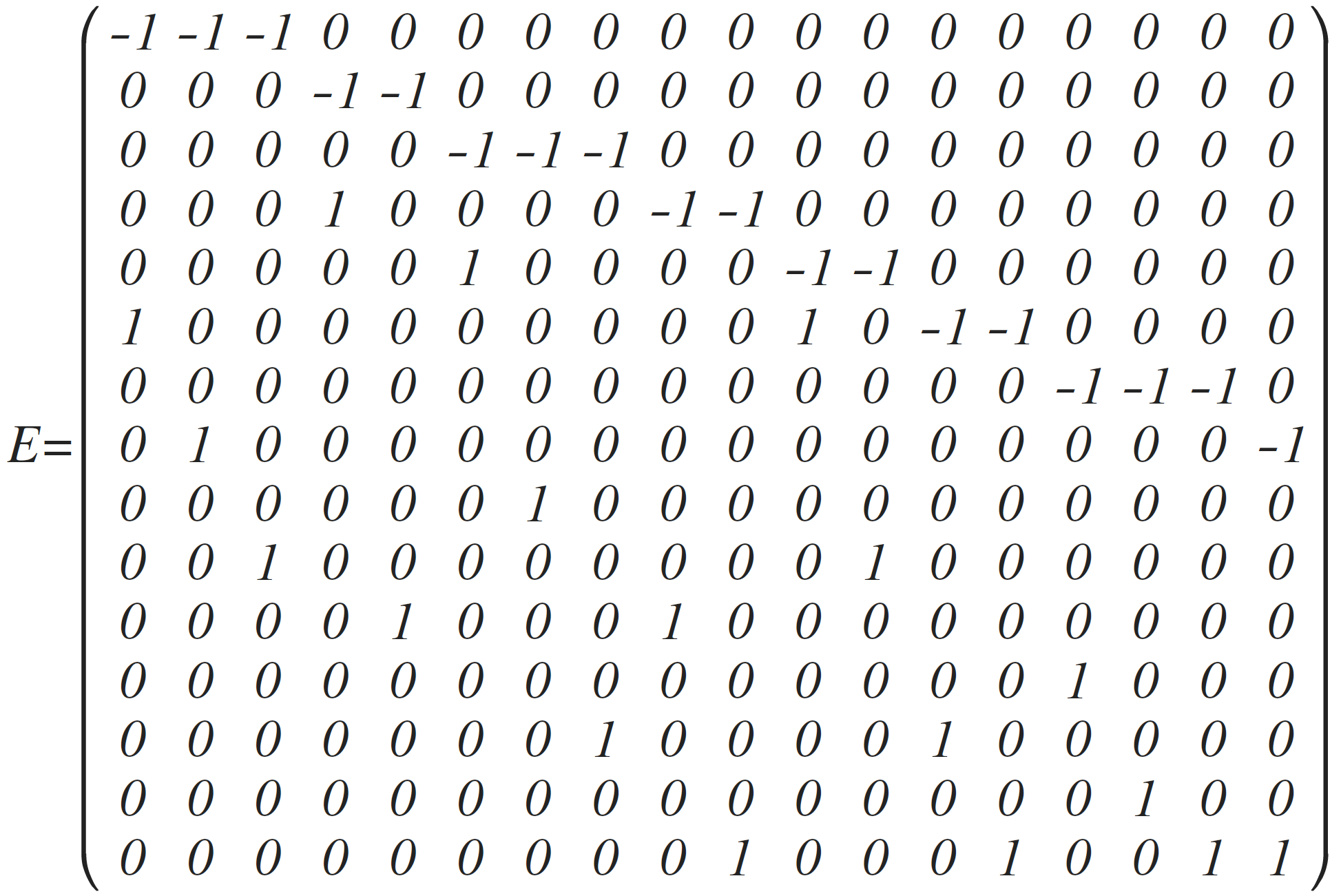}
  \caption{Incidence matrix $E$.}  \label{fig:E}
\end{figure}

\begin{figure*}[t]
  \centering
     \begin{subfigure}[b]{0.45\textwidth}
         \centering
  \includegraphics[scale=0.7]{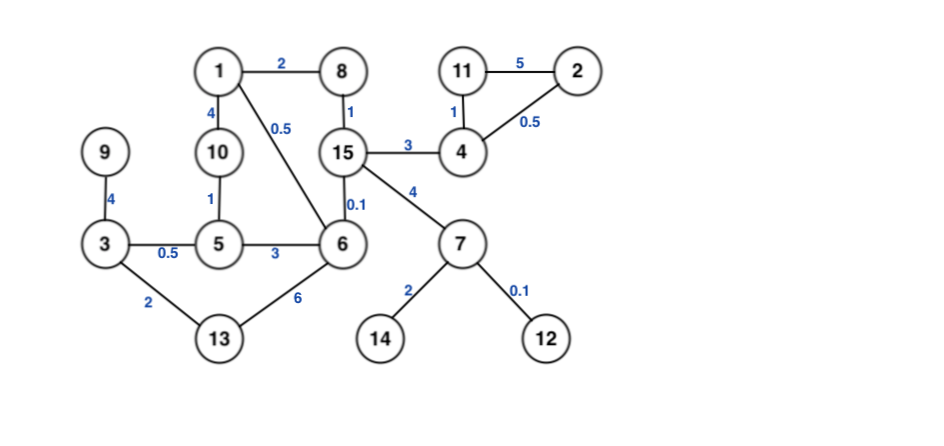}
    \caption{Connected network with 15 nodes. The numbers above arrows denote the weights of the edge between any two agents.}
          \label{fig:g1}
     \end{subfigure}
     \hfill
     \begin{subfigure}[b]{0.45\textwidth}
         \centering
          \includegraphics[scale=0.35]{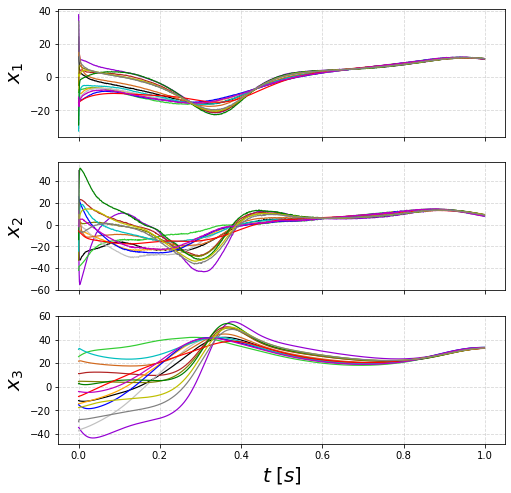}
         \caption{Evolution of the states $(x_1,x_2,x_3)$ of agents $a_i, \, i=1,\dots, N$.}
    \label{fig:sim1}
     \end{subfigure}
        \caption{Synchronization of Lorentz's oscillators}
        \label{fig:lorenz}
\end{figure*}

Building up on similar robustness properties, 
we illustrate in this section that 
we can learn the feedback $\alpha$
of Theorem~\ref{the1} and yet achieving
asymptotic synchronization of a complex undirected network.
To this end,  we use the Deep Neural Networks (DNN) based algorithm for metric estimation presented in \cite[Section VII-B]{giaccagli2024synchronization} and we apply to 
a network of $N$ identical (chaotic) Lorenz attractors. Each agent $i\in\{1,\ldots, N\}$ is described by the following dynamics
$$
\left\{
\begin{array}{l}
\dot{x}_{i,1}  =  a(x_{i,2}-x_{i,1})+u_i \\
\dot{x}_{i,1}  =  x_{i,1}(b-x_{i,3})-x_{i,2}+(2+\sin(x_{i,1}))u_i \\
\dot{x}_{i,3}  =  x_{i,1}x_{i,2}-cx_{i,3}
\end{array}
\right.
$$
where $x_i=(x_{i,1},x_{i,2},x_{i,3})\in \RR^3$ is the state of agent $i$, $u_i\in \RR$ the control input, 
and the parameters $a,b,c$ are slected as
$a=10$, $b=\frac{8}{3}$, and $c=28$ in order to guarantee 
chaotic behaviour.
Note that with such a choice of control vector field $g$, the system 
is not feedback linearizable and conditions
of Assumption~\ref{ass:riem} and
~\ref{ass:killing} cannot be easily verified 
due to the presence of both 
polynomial and trigonometric functions.
A DDN is used to learn an approximated solution on a sufficiently large region following \cite[Section VI]{giaccagli2024synchronization}.
In the illustration, we considered an undirected and connected network of $N=15$ agents.
The network is illustrated in Figure \ref{fig:g1}, where the numbers over each arcs are the weight of each connections.
For this example we have the incidence matrix depicted in 
Figure~\ref{fig:E}
with a  matrix of the weight of each nodes  $W = \diag(0.5,2,4,0.5,5,0.5,4,2,1,3,3,1,6,0.1,0.1,2,4,1)$.
Figure~\ref{fig:sim1} show that if the initial conditions are not too far, all the agents of the network synchronize among them with the distributed control law \eqref{ui} derived in  Theorem~\ref{the1}.

\color{black}

\section{Conclusion} \label{sec:ccl}
Based on a new technical result
on the properties of the Edge Laplacian 
matrix of a graph, we have been able 
to prove  that the 
distributed feedback
design 
proposed in 
\cite{giaccagli2021sufficient,giaccagli2024synchronization}
is able to solve the problem
of global state synchronization 
of a network of 
input-affine nonlinear agents
connected via a 
weighted undirected graph.
The proof is based on the construction 
of a new Lyapunov function 
which is formally equivalent to the 
Edge's Riemannian energy.
It is interesting to note that the derivative of the proposed Lyapunov function is always negative regardless the 
connectivity properties of the graph.
However, it defines a suitable Lyapunov function 
for the network synchronization only when the graph is connected.

We remark that similarly to \cite{giaccagli2024synchronization},
the Killing and integrability conditions
of Assumption \ref{ass:killing}
can be relaxed at the price of 
guaranteeing only 
semi-global synchronization.
In turns, such a relaxation  allows the use of  an approximation
of the feedback law
$\alpha$ which is generically hard to obtain in analytic form.
Such an approximation can be computed
via numerical methods such as 
neural networks  (e.g. 
\cite{tsukamoto2021contraction,giaccagli2024synchronization}) while guaranteeing asymptotic synchronization. 

Finally, we conjecture that the proposed
feedback law can apply also to the more general case
of weighted
directed graphs in the analysis.
As a future work, we would like to certify this conjecture by a 
suitable Lyapunov analysis (which may include also 
the use of colour theory as in 
\cite{pavlov2022nonlinear}).

\bibliography{biblio.bib}

\end{document}